\documentclass[11pt]{article}

\usepackage[numbers]{natbib}
\usepackage{amsfonts}
\usepackage{amsmath}
\usepackage{amsthm}
\usepackage{amssymb}
\usepackage{setspace}
\theoremstyle{plain}
\newtheorem{thm}{Theorem}
\newtheorem{defin}{Definition}
\newtheorem{prop}{Proposition}

\begin{document}

\title{On regularities of mass random phenomena}%

\author{Victor I. Ivanenko \footnote{victor.ivanenko.1@gmail.com} \footnote{Department of mathematical modeling of economical systems, Kyiv Polytechnical Institute, 57, Prospekt Peremogy, Kyiv 03056, Ukraine,}, \hspace{10pt} \fbox{Valery A. Labkovsky}} 

\date{\vspace{-5ex}}

\maketitle

\begin{abstract}
This paper contains an answer to the question of existence of regularities of the so called \textit{random in a broad sense} mass phenomena, asked by A. N. Kolmogorov in \cite{Kolmogorov}. It turns out that some family of finitely-additive probabilities is the statistical regularity of any such phenomenon. If the mass phenomenon is stochastic, then this family degenerates into a single probability measure. The paper provides  definitions, the formulation and the proof of the theorem of existence of statistical regularities, as well as the examples of their application.
\\
\textbf{Key words} \textit{Sequence, Net, Nonstochastic randomness, Statistical regularity, Families of probability distributions}
\end{abstract}

\section{Introduction}
The interest in studying the properties of \textit{mass random phenomena} (MRP) is not new. So in \cite{Jarvik, Borel} the authors pointed out the difficulties arising in the process of modeling of the social MRP. In particular, in \cite{Borel} one reads:"Some contemporary theoreticians ... think that probability could be defined as frequency for a very large number of trials. If for a very big number of trials this frequency does not tend to a limit, but fluctuates more or less between different limits, one needs to affirm that probability $p$ does not remain constant and changes in the process of trials. This concerns, for example, human mortality rate in the course of centuries, since the progress of medicine and hygiene leads to the increase of life duration." The problem of revealing of the regularities of MRP becomes more and more important, especially in relation to the instability of financial markets and other economic objects \cite{Taleb,Mandelbrot, Munier}, that makes forecasting in this area very unreliable. 

In \cite{Khinchin1} the question was risen whether the MRP posses the properties that are necessary in order to apply the probability theory to their description?

In \cite{Kolmogorov} we find the following remark by A.N. Kolmogorov: "Speaking of randomness in the ordinary sense of this word, we mean those phenomena in which we do not find regularities allowing us to predict their behavior. Generally speaking, there are no reasons to assume that random in this sense phenomena are subject to some probabilistic laws. Hence, it is necessary to distinguish between randomness in this \textit{broad sense} and \textit{stochastic} randomness (which is the subject of probability theory)". 

However, what do the words "do not find regularities allowing us to predict their behavior" mean? Hardly these words should be understood in the sense that such regularities do not exist at all. More likely, these words point out to the problem of finding statistical regularities of mass \textit{random in a broad sens phenomena} (MRBSP), that is the regularities of asymptotic behavior of different average values that characterize these phenomena. It can be frequencies of hitting in given subsets, arithmetic averages of some functionals, and so on. Recall that MRP are called \textit{statistically stable} or \textit{stochastic}, if with the increase of the number of "trials" all these averages tend to limits (and if some other conditions are verified as well, see details in \cite{Kolmogorov}). Unlike this, it is natural to consider as MRBSP those MRP, whose behavior interests us only to within their  statistical regularities. In other words, this definition combines in MRBSP stochastic as well as \textit{nonstochastic}  random phenomena. \footnote{The term "nonstochastic" appeared in \cite{Vyugin} in the context of Kolmogorov's complexity, meaning "more complex than stochastic". In this paper the meaning of this term is  "more random than stochastic".} 

The theorem of existence of statistical regularities of these phenomena in the form of families of probability distributions, and their significance to decision theory, constitute the content of this paper.


\section{Theorem of existence of statistical \\regularities}

An ordinary sequence is the simplest mathematical model of a mass phenomenon. In order to construct, on the basis of a sequence, a model of a random phenomenon, it is necessary to identify sequences that have identical  statistical properties.

\begin{defin}
Let $X$ be an arbitrary set. Two sequences $\overline{x}^{(1)}$ and $\overline{x}^{(2)}$ of elements of the set $X$ are called  statistically equivalent ($S$-equivalent) if and only if for any natural number $m$ and any bounded mapping $\gamma \in (X \rightarrow \mathbb{R}^m)$ the set of limit points of the sequence 
\begin{eqnarray*}
\left\{\overline{y}_n^{(k)}; n \in \mathbb{N}\right\}, \overline{y}_n^{(k)} =\frac{1}{n}\sum_{i=1}^{n}\gamma(\overline{x}_{i}^{(k)})
\end{eqnarray*}
does not depend on $k \in \left\{1,2\right\}$.
\end{defin}

The class of $S$-equivalence of the sequence $\overline{x} \in X^{\mathbb{N}}$ will be denoted as $S(\overline{x})$. Our nearest goal is to find the invariant of the relation of $S$-equivalence. Introduce several notions.

Let $M$ be a Banach space of bounded real functions, defined on the set $X$, $M^{*}$ be the dual space of the space $M$, and $\tau$ - is a weak-$*$ topology in $M^{*}$. Let, further, $PF(X)$ be the subspace of the topological space $(M^{*},\tau)$ defined by the formula 
\begin{eqnarray*}
PF(X)=\left\{p \in M^{*}: p(f)\geq 0 \text{ if } f \geq 0, p(\mathbf{1}_X)=1 \right\},
\end{eqnarray*}
where $\mathbf{1}_A(\cdot)$ is the characteristic function of the set $A$.\footnote{In what follows, instead of $p(\mathbf{1}_A)$ we shall often write $p(A)$,  identifying, by the same token, the \textbf{elements} of the set $PF(X)$ with the finitely additive and normed measures on $2^X$. Obviously, $p(f)$ in this case is simply the integral $p(f)=\int{f(x)p(dx)}$, defined naturally due to boundedness of function $f$.}

Associate to an arbitrary sequence $\overline{x}=\left\{\overline{x}_n; n \in \mathbb{N}\right\} \in X^{\mathbb{N}}$ the  sequence of \textbf{elements} from $PF(X)$ defined as
\begin{eqnarray*}
\left\{\overline{p}_{\overline{x}}^{(n)}(\cdot); n \in \mathbb{N}\right\}, \overline{p}_{\overline{x}}^{(n)}(A)=\frac{1}{n}\sum_{i=1}^{n}\mathbf{1}_A(\overline{x}_i), \forall A \subseteq X.
\end{eqnarray*}
\textbf{Note that the \textbf{elements} $\overline{p}_{\overline{x}}^{(n)}(A)$ from $PF(X)$ has the meaning of frequencies of hitting of the elements of $\overline{x}$ in the set $A$}.\footnote{See also \cite{Mises,Kolmogorov1}} Due to compactness of the space $PF(X)$ (as of a bounded closed set in $(M^*,\tau)$), the sequence $\left\{\overline{p}_{\overline{x}}^{(n)}(\cdot); n \in \mathbb{N}\right\}$ will have a non-empty closed set of limit points, which we denote as $P_{\overline{x}}$ and call \textit{the regularity} of this sequence. Therefore introduce the following definition.

\begin{defin}
Any non-empty closed subset of the space $PF(X)$ is called \textbf{a regularity} on $X$. Denote the set of all regularities on $X$ as $\mathbb{P}(X)$  and associate to any sequence $\overline{x} \in X^{\mathbb{N}}$ its regularity $P_{\overline{x}}$. Finally, for $m \in \mathbb{N}$, $\gamma=(\gamma_1,\gamma_2,\dots,\gamma_m)\in (X \rightarrow \mathbb{R}^m)$ and $P \in \mathbb{P}(X)$, the symbol $P(\gamma)$ denotes the set 
\begin{center}
$\left\{ (r_1, r_2,\dots,r_m) \in \mathbb{R}^m: \exists p \in P , r_i=p(\gamma_i), \forall i \in \overline{1,m} \right\}$, 
\end{center}
and, in particular, $p(\gamma)=(p(\gamma_1),p(\gamma_2),\dots,p(\gamma_m))$ for $p \in PF(X)$.
\end{defin}

Consider the following proposition.
\begin{prop}
The mapping $\overline{x} \mapsto P_{\overline{x}}$ is the invariant of the relation of $S$-equivalence on $X^{\mathbb{N}}$.
\end{prop}

This statement will be proved below in a more general form. So far, however, let us agree to call the classes of $S$-equivalence of sequences \textit{the simplest random phenomena}, and their regularities - \textit{statistical regularities} of the corresponding phenomena. Any sequence $\overline{x} \in X^{\mathbb{N}}$ is considered as a realization of a simplest random phenomenon $S(\overline{x})$.

Connection of the notions introduced above with the probabilistic notions follows directly from the enforced law of large numbers.
\begin{prop}
Let $X$ be a finite set, $\mu$ - a probability distribution on $X$, and $\overline{\xi}=\left\{ \overline{\xi}_n; n \in \mathbb{N} \right\}$ - a sequence of independent (in the usual sense) random elements, taking values in $X$ with distribution $\mu$. Then with probability $1$ the sequence $\overline{x}$ of the values of the sequence $\overline{\xi}$ will be a realization of the simplest random phenomenon $S(\overline{x})$ with statistical regularity $P_{\overline{x}}=\left\{ \mu \right\}$, i.e.  consisting of the single distribution $\mu$.
\end{prop}

However, when the set $X$ is infinite everything becomes considerably more difficult. In this case, the capabilities of sequences, generally speaking, are insufficient in order to guarantee that the frequencies of hitting in all measurable sets would tend to their limits simultaneously. Moreover, it is easy to see that the regularities of sequences, since they are concentrated only on a countable subset of the set $X$, constitute only a small part of the set of all regularities on $X$. This seems to reflect the fact that sequences constitute only a small part of all mass phenomena. A more general notion of \textit{sampling net} is, as we shall see further, already sufficient for our goals. 

\begin{defin}
A \textbf{sampling net} (s.n.) in $X$ is any net $\varphi=\left\{\varphi_{\lambda}, \lambda \in \Lambda, \geq  \right\}$ taking values in the sampling space 
\begin{eqnarray*}
X^{\infty}=\bigcup_{n=1}^{\infty}X^n, \hspace{0.25cm} X^n=\underbrace{X \times \cdots \times X}_{n}.
\end{eqnarray*}
Moreover, if $\lambda \in \Lambda$,$ \varphi_{\lambda} \in X^n$ then we denote $n=n_{\lambda}, \varphi_{\lambda}=(\varphi_{\lambda 1}, \varphi_{\lambda 2},\dots, \varphi_{\lambda n_{\lambda}}) $ and associate to this $\lambda$ the measure $p_{\varphi}^{(\lambda)} \in PF(X)$ defined as
\begin{eqnarray*}
p_{\varphi}^{(\lambda)}(A)=\frac{1}{n_{\lambda}}\sum_{i=1}^{n_{\lambda}}\mathbf{1}_A(\varphi_{\lambda i}),\quad A \subseteq X.
\end{eqnarray*}
The set $P_{\varphi}$ of limit points of the net $p_{\varphi}=\left\{p_{\varphi}^{\lambda}, \lambda \in \Lambda, \geq  \right\}$ will be called \textbf{the regularity} of the s.n. $\varphi$. The class of all s.n. in $X$ will be denoted as $\Phi(X)$.
\end{defin}

Extend now the relation of $S$-equivalence on the whole $\Phi(X)$.

\begin{defin}
Sampling nets $\varphi^{(k)} \in \Phi(X)$,$k=1,2$ are considered as $S$-equivalent if and only if for any $m \in \mathbb{N}$ and any bounded mapping $\gamma \in (X \rightarrow \mathbb{R}^m)$ the set of limit points of the net of averages 
\begin{eqnarray}\label{s_directedness}
\left\{ y_{\lambda}^{(k)}, \lambda \in \Lambda, \geq \right\}, \quad y_{\lambda}^{(k)}=\frac{1}{n_{\lambda}}\sum_{i=1}^{n_{\lambda}}{\gamma(\varphi^{(k)}_{\lambda i})}
\end{eqnarray}
does not depend on $k \in \left\{ 1,2 \right\}$.
\end{defin}

We can now formulate the main theorem in the following way.

\begin{thm}\label{rbs_theorem}
\begin{enumerate}
\item [(i)]For any s.n. $\varphi \in \Phi(X)$, any $m \in \mathbb{N}$ and any bounded mapping $\gamma \in (X \rightarrow \mathbb{R}^m)$, the set of limit points of the net (\ref{s_directedness}) can be written as $P_{\varphi}(\gamma)$.
\item [(ii)]The mapping $\varphi \mapsto P_{\varphi}$, defined on $\Phi(X)$, is the invariant of the relation of $S$-equivalence.
\item [(iii)] This mapping is a mapping on the whole set $\mathbb{P(X)}$, i.e. the set $\Phi(X)/S$ of classes of $S$-equivalence and the set $\mathbb{P}(X)$ of regularities are put by this mapping into one-to-one correspondence.
\end{enumerate}
\end{thm}

This theorem justifies the following definition.

\begin{defin}\label{rbs}
Any class of $S$-equivalence of sampling nets in $X$ is called random in a broad sense mass phenomenon in $X$. The regularity $P_{\varphi}$ is called the statistical regularity of the random phenomenon $S(\varphi)$. Any s.n. $\varphi^{'} \in S(\varphi) $ is called a realization of the random phenomenon $S(\varphi)$. The random phenomenon, having statistical regularity $P$, is called $\mu$-stochastic if and only if there exists a non-trivial $\sigma$-algebra $\mathcal{A} \subseteq 2^X$, on which $\mu$ is a $\sigma$-additive probability, and $p(A)=\mu(A)$ for all $p \in P, A \in \mathcal{A}$.
\end{defin}

\section{The proof}

Denote the set of limit points of an arbitrary net $g=\left\{g_\alpha, \alpha \in A, \succcurlyeq \right\}$  with values in $X$ as $LIM(g)$ or $LIM\left\{g_\alpha, \alpha \in A, \succcurlyeq \right\}$. Denote the set of bounded mappings from $X$ into $\mathbb{R}^m$ as $M^m$. We need to establish the three following facts:
\begin{itemize}
\item[$(i)$]\hspace{0.25cm}The relation
\begin{equation*}
LIM\left\{y_\lambda, \lambda \in \Lambda, \succcurlyeq \right\}=P_{\varphi}(\gamma), 
\end{equation*}
where
\begin{equation*}
y_\lambda=\frac{1}{n_\lambda} \sum_{i=1}^{n_\lambda}\gamma(\varphi_{\lambda i}),
\end{equation*}
is true for all $m \in \mathbb{N}$, $\gamma \in M^m$, $\varphi \in \Phi(X)$.
\item[$(ii)$]\hspace{0.25cm} If $P_1, P_2 \in \mathbb{P}(X), P_1 \neq P_2$, then there exist such $m \in \mathbb{N}$ and such $\gamma \in M^m$, that $P_1(\gamma) \neq P_2(\gamma)$.
\item[$(iii)$]\hspace{0.25cm} For any regularity $P \in \mathbb{P}(X)$ there exist such s.d. $\varphi \in \Phi(X)$, that $P = P_{\varphi}$.
\end{itemize}
Begin with the proof of the proposition $(i)$. Let $r \in LIM(y)$, where $y=\left\{y_\lambda, \lambda \in \Lambda, \succcurlyeq \right\}$. Then there exists a subnet of the net $y$ converging to $r$, i.e. there exists (see \cite{Kelley}) a directed set $(A, \succcurlyeq)$ and a function $f: A \rightarrow \Lambda $ such that the net $\overline{y}=y \circ f$ converges to $r$, and, in addition, for any $\lambda \in \Lambda$ there exists such $\alpha_1 \in A$ that $f(\alpha) \succcurlyeq \lambda$ for all $\alpha \succcurlyeq \alpha_1$.

Consider now the net of measures $\overline{p}_{\varphi}=p_{\varphi} \circ f$, where $ p_{\varphi}=\left\{p^{(\lambda)}_{\varphi}, \lambda \in \Lambda, \succcurlyeq \right\}$. By virtue of compactness of the space $(PF(X), \tau)$ this net has at least one limit point. Denote it as $p_0$ and consider a subnet  $\overline{\overline{p}}_{\varphi}$ of the net $\overline{p}_{\varphi}$, converging to $p_0$. Let it be $\overline{\overline{p}}_{\varphi}=\overline{p}_{\varphi} \circ g = p_{\varphi} \circ f \circ g, g : B \rightarrow A$. Then the net $\overline{\overline{y}}= y \circ f \circ g$, on the one hand, converges to $r$, and, on the other hand, $\overline{\overline{y}}_{\beta}=\overline{\overline{p}}_{\varphi}^{(\beta)}(\gamma), \beta \in B$, so that 
\begin{eqnarray*}
r=\lim_{\beta}\overline{\overline{p}}_{\varphi}^{(\beta)}(\gamma)=p_0(\gamma) \in P_{\varphi}(\gamma).
\end{eqnarray*}
By the same token, it is proved that $LIM(y) \subseteq P_{\varphi}(\gamma)$.

Conversely, if $p_0 \in P_{\varphi}, r=p_0(\gamma)$, then there exists a subnet $\widetilde{p}_{\varphi}=\left\{\widetilde{p}_{\varphi}^{\alpha},\alpha \in A, \succcurlyeq  \right\}$ of the net $p_{\varphi}$, converging to $p_0$. But in this case $\lim_{\alpha}\widetilde{p}_{\varphi}^{(\alpha)}(\gamma_i)=p_0(\gamma_i)$ for all $i \in \overline{1,m}$. It means that $\lim_{\alpha}\widetilde{p}_{\varphi}^{(\alpha)}(\gamma)=p_0(\gamma)$. And, since $\widetilde{p}_{\varphi}^{(\alpha)}(\gamma)=y_\lambda$ for $\lambda=f(\alpha)$, this proves $(i)$.

In order to prove $(ii)$ assume that there exists $p_1 \in P_1 \setminus P_2$. Since the set $P_2$ is closed, there exists a vicinity of the point $p_1$ that does not cross with $P_2$ and it means that there exist such $\epsilon > 0, \gamma_1,\gamma_2,\dots,\gamma_m \in M$ that
\begin{eqnarray*}
\forall p_2 \in P_2, \exists i \in \overline{1,m}, \left| p_1(\gamma_i) -  p_2(\gamma_i) \right| > \epsilon.
\end{eqnarray*}
So that if $\gamma=(\gamma_1,\gamma_2,\dots,\gamma_m)$, then $p_1(\gamma) \notin P_2(\gamma)$.

The complete proof of $(iii)$ can be found in \cite{IvanenkoLabkovsky1,IvanenkoMunier1,Ivanenko}. Here we shall outline the main ideas of the proof. Let $Q$ be the set of all such measures $q \in PF(X)$ that each one of them is concentrated on a finite set $X_q \subseteq X$, and in addition all numbers $q(x), x \in X_q$ are rational. One can show that the set $Q$ is everywhere dense in $(PF(X),\tau)$.

Now, to an arbitrary regularity $P \in \mathbb{P}(X)$ we put into correspondence the directed set $(\Lambda, \succcurlyeq)$ such that 
\begin{eqnarray*}
\Lambda=\mathbb{R}^{+} \times M^{\infty} \times P, \hspace{6pt} \mathbb{R}^{+}=]0,\infty[, \hspace{6pt} M^{\infty}=\bigcup_{m=1}^{\infty}M^m,M^m=\underbrace{M \times \dots \times M}_{m}
\end{eqnarray*}
and the relation $(\succcurlyeq)$ is given by the formula
\begin{eqnarray*}
(\epsilon_1,\gamma_{11},\gamma_{12},\dots,\gamma_{1n_{1}}, p_1) \succcurlyeq (\epsilon_2,\gamma_{21},\gamma_{22},\dots,\gamma_{2n_{2}}, p_2) \Leftrightarrow \\ (\epsilon_1 \leq \epsilon_2, \left\{ \gamma_{11},\gamma_{12},\dots,\gamma_{1n_{1}} \right\} \supseteq \left\{ \gamma_{21},\gamma_{22},\dots,\gamma_{2n_{2}} \right\}),
\end{eqnarray*}
where no condition is imposed on $p_1$ and $p_2$. 

Finally, to any $\lambda = (\epsilon,\gamma_{1},\gamma_{2},\dots,\gamma_{m}, p) \in \Lambda$ we put into correspondence some
\begin{eqnarray*}
q_\lambda \in Q \bigcap \left\{ p^{'} \in PF(X): \forall i \in \overline{1,m}, \left| p(\gamma_i)  - p^{'}(\gamma_i) \right| < \epsilon \right\}.
\end{eqnarray*}

It is proven further that with any $\lambda \in \Lambda$ one can associate simultaneously a sequence of points $x_1^{(\lambda)}, x_2^{(\lambda)},\dots,x_{n_\lambda}^{(\lambda)} \in X_q$ satisfying the condition
\begin{eqnarray*}
q_\lambda (A)=\frac{1}{n_\lambda}\sum_{i=1}^{n_\lambda} \mathbf{1}_A(x_i^{(\lambda)}), \quad \forall A \subseteq X.
\end{eqnarray*}
It remains to chose $\varphi_\lambda=(x_1^{(\lambda)}, x_2^{(\lambda)},\dots,x_{n_\lambda}^{(\lambda)})$ and we obtain a s.n. $\varphi: \lambda \mapsto \varphi_\lambda$ that has the regularity $P_\varphi=P$.

\section{Applications in decision theory}

Statistical regularities of the general form find their application in decision theory \cite{IvanenkoLabkovsky0, IvanenkoLabkovsky11, Gilboa, IvanenkoLabkovsky1, IvanenkoMunier1, Ivanenko, Mikhalevich} and its applications \cite{IvanenkoMunier2}.

Considering decision problems, assume that we need to make a decision $u$ from the set $U$ of possible decisions, knowing that the result of making a decision depends on some uncontrolled parameter $\theta$ from the set $\Theta$ of possible values of this parameter and is described by the bounded real loss function  $L: \Theta \times U \rightarrow \mathbb{R}$. 
If nothing is known about the behavior of the parameter, then we cannot, strictly speaking, exclude that scenario, where the value of $\theta$ is chosen in the worst possible for us way. In this case, the quality of decision  $u$ is evaluated by means of the loss function 
\begin{eqnarray*}
L_1^* (u)=\sup_{\theta \in \Theta}L(\theta,u), \hspace{0.25cm} u \in U,                            
\end{eqnarray*}
a so called "minmax" criterion.

If it is known, that parameter $\theta$ is stochastic with the given distribution $\mu$, then, trying to minimize the average losses, one makes use of the Bayes criterion
\begin{eqnarray*}
L_2^* (u)=\int L(\theta,u) \mu(d\theta),\hspace{0.25cm} u \in U.                            
\end{eqnarray*}

Suppose now that parameter $\theta$ is random in a broad sense with the statistical regularity $P \in \mathbb{P}(\Theta)$. Let us show that in this case it is natural to chose the criterion in the form of 
\begin{eqnarray}\label{criterion}
L_3^* (u)=\sup_{p \in P}\int L(\theta,u)p(d\theta), \hspace{0.25cm} u \in U,                            
\end{eqnarray}

Indeed, let $r_1<L_3^*(u)<r_2$. The following statement is straightforward
\begin{prop}\label{proposition3}
Let $\left\{ \varphi_{\lambda}, \lambda \in \Lambda, \succcurlyeq \right\}$ - be a sampling net in $\Theta$ with the regularity $P$. Then for any $\lambda_1 \in \Lambda$ there is such $\lambda \succcurlyeq \lambda_1$ that
\begin{eqnarray*}
\frac{1}{n_{\lambda}}\sum_{i=1}^{n_{\lambda}}{L(\varphi_{\lambda i}, u)} > r_1
\end{eqnarray*}

and, at the same time, there is such $\lambda_2$, that for all $\lambda \succcurlyeq \lambda_2$ there will be 
\begin{eqnarray*}
\frac{1}{n_{\lambda}}\sum_{i=1}^{n_{\lambda}}{L(\varphi_{\lambda i}, u)} < r_2.
\end{eqnarray*}

\end{prop}
In other words, $L_3^* (u)$ - is that natural border, that separates the average losses, that can happen for a given $u$ for an arbitrary "large" $\lambda$, from those average losses that are not "dangerous" to us, when $\lambda$ is sufficiently "large".

It is easy to see that $L_3^* (u)$ becomes $L_1^* (u)$, when $P=PF(\Theta)$ (strictly nothing is known about $\theta$, save the set $\Theta$ where it takes values), and that it becomes $L_2^* (u)$, when $P={\mu}$ is stochastic regularity and function $L(\cdot,u)$ is measurable relatively to the corresponding $\sigma$- algebra.

The inverse result appears as somewhat surprising. It turns out that if one subordinates a criterion choice rule to some natural conditions of consistency with the triplet  $(\Theta,U,L)$, then any rule, satisfying these conditions, leads to the criterion of the form (\ref{criterion}), where $P$ - is some (not known beforehand) regularity on $\Theta$. In particular, this result justifies the heuristic definition of random in a road sense phenomena introduced above. Therefore, one can conclude that regularity on $\Theta$ is, in a certain sense, the most general form of information about the behavior of $\theta$. One can find details in \cite{IvanenkoLabkovsky0, IvanenkoLabkovsky11,IvanenkoLabkovsky1, IvanenkoMunier1, Ivanenko}.

\section{Concluding remarks}

Today there are several approaches to modeling of the MRP. So, there is the algorithmic approach to randomness \cite{Zvonkin, Vyugin} as well as the game-theoretic approach to randomness in finance \cite{Vovk}. An alternative approach was studied in \cite{IvanenkoLabkovsky3}, where sequences were constructed  only with the requirement of the so called $\Gamma_{l}$- independence. In mathematical finance diverse extensions of stochastic models have been  popular \cite{Heston, Avellaneda, Calvet}. 

Families of probability distributions appear in literature more and more often. So, they were considered in game theory \cite{Shapley} in order to study non-additive set functions. In the so called subjective decision theory these families appear as consequence of the axioms of rational choice  \cite{IvanenkoLabkovsky0,IvanenkoLabkovsky11, Gilboa}, where, similarily to robust statistics \cite{Huber}, they were interpreted as families of a priori distributions.  Families of probability distributions attract all the more attention in statistical data analysis  \cite{Fine}\footnote{I am thankful to professor Vladimir Vovk who made me familiar with the works of professor Terrence Fine and, in particular, with this paper.}.

It turns out that specifically families of probability distributions are necessary for the description of statistical (frequentist) regularities of a rather wide class of MRBSP. The theorem of existence of such statistical regularities was published earlier in somewhat different form \cite{IvanenkoLabkovsky2, IvanenkoLabkovsky1}.

\end{document}